\documentclass[11pt]{article} 
\usepackage{latexsym,amssymb}

\font\smc=cmcsc10
\font\tenbbb=msbm10
\font\sevenbbb=msbm7
\font\fivebbb=msbm5
\newfam\msbfam
\textfont\msbfam=\tenbbb \scriptfont\msbfam=\sevenbbb
\scriptscriptfont\msbfam=\fivebbb
\def\bbb{\fam \msbfam \tenbbb}
\font\tenfrak=eufm10
\font\sevenfrak=eufm7
\font\fivefrak=eufm5
\newfam\euffam
\textfont\euffam=\tenfrak \scriptfont\euffam=\sevenfrak
\scriptscriptfont\euffam=\fivefrak

\def\br{\hfill\break}
\def\3#1#2#3{_{#1=#2}^{#3}}	

\def\N{{\bbb N}}

\def\3#1#2#3{_{#1=#2}^{#3}}	
\def\implies{\; \Longrightarrow \;}

\def\mdef#1#2{\newcommand#1{\ensuremath{#2}}}
\def\kdots{\ifmmode,\ldots,\else,~$\ldots\,$, \fi}
\def\row#1#2{\ifmmode #1_1\kdots #1_{#2}\else
 $#1_1$\kdots$#1_{#2}$\fi}
\def\orow#1#2{\ifmmode #1_0\kdots #1_{#2}\else
 $#1_0$\kdots$#1_{#2}$\fi}

\def\sqr#1#2{{\vcenter{\hrule height.#2pt
\hbox{\vrule width.#2pt height#1pt \kern#1pt
\vrule width.#2pt}
\hrule height.#2pt}}}
\def\square{\ifmmode\msquare\else\hfill$\sqr55$\medskip\fi}
\def\msquare{\eqno\sqr55}
\long\def\proclaim #1. #2\par{\medbreak
 \noindent {\smc #1.\enspace }{\sl #2}\par
 \ifdim \lastskip <\medskipamount \removelastskip
 \penalty 55\medskip \fi}
\def\proc #1. {\medbreak\noindent {\smc #1.\enspace }}
\long\def\clause#1#2{\par\smallbreak\hangafter=1\hangindent
 20pt\noindent{\hbox to20pt{{#1\hfill}}{#2}}\hfill\par
 \ifdim \lastskip <\smallskipamount \removelastskip
 \penalty 15\smallskip \fi}

\newcommand{\pureca}[1]{\left|{#1}\right|}
\def\ca#1 {\ensuremath{\pureca{#1}}}
\def\comp#1{\mathop{\rm comp}({#1})}

\mdef\0{\widehat0}
\mdef\1{\widehat1}
\def\A{\ensuremath{{\cal A}}}

\def\P{{\bbb P}}
\def\rank{\mathop{\rm rank}\nolimits}
\newcommand{\2}[2]{\ensuremath{\mathopen[#1,#2\mathclose]}}
\let\length=\ell

\begin{document}

\title{Binomial posets with non-isomorphic intervals}

\author{{\sc J\"orgen Backelin}}

\date{}

\maketitle

\begin{abstract}
The confluent
binomial posets with the atomic function
$A(n) = \max(1,2^{n-2})$ are classified. In particular, it is
shown that in general there are
non-isomorphic intervals of the same length.
\end{abstract}



\renewcommand{\c}[1]{\ensuremath{C_{#1}}}
 \section{Introduction}

 Recall that (following \cite{Stanley_a}) a binomial poset
$(P,\le)$ is a graded locally finite poset, where the number
of maximal chains in an interval only depends on the length of
that interval. In this article, we always consider directed
infinite downwards bounded binomial posets. In other word,
any poset $P$ will have a minimal element \0; for any $x,y
\in P$, there is a $z \in P$ with $x \le z$ and $y \le z$;
for any (by definition non-empty) interval $\2xy \subset P$;
each interval is a finite set;
all inclusion maximal chains in \2xy have the same length
$\length(\2xy)$; and the number of such chains for a given interval
\2xy is only depending on $P$ and $\length(\2xy)$. In fact, we
mostly consider \emph{(strongly) confluent} such posets, where
in addition there is an infinite chain $x_1 < x_2 < x_3 < \ldots$
in $P$, such that $P = \bigcup_i \2\0{x_i}$. (The reader may verify
that a downward bounded directed locally finite poset is strongly
confluent if and only if it is countable as a set.)

 Call \2xy an
$n$-interval, if $\length(\2xy) = n$; and let $B(P,n)$ (or $B(n)$,
if there is no ambiguity) denote the number of maximal chains
in an $n$-interval of $P$. The numbers $B(P,0)=1$, $B(P,1)=1$,
$B(P,2)$, $B(P,3)$\kdots behave as a kind of generalised
factorial functions, in the following senses: For any
$n$-interval \2xy, the number of elements therein that cover
$x$ is an integer $A(P,n)$ (or $A(n)$), only depending on $P$
and $n$; $B(P,n) = \prod\limits\3i1nA(P,n)$; and for any
integers $j$ and $n$, such that $0 < j < n$, the number
$B(P,n) \big/ \bigl( B(P,j)\,B(P,n-j) \bigr)$ counts the
number of 2-chains $x<z<y$ in \2xy, such that $z$ has local
rank $j$ (i.e., such that $\length(\2xz)
= j$), and thus in particular is a positive integer, a
`generalised binomial coefficient'.

 In the sequel, we ignore $B(0)=1$ and $A(0)=0$, but consider
posets $X$ with a given fixed sequence $B(1)$, $B(2)$\kdots and
thus a fixed atomic sequence
$A(1) = a_1 = 1$, $A(2) = a_2$, {\sl et cetera}. The
{\sl atomic sequence\/} $\A := (a_i)\3i1\infty$ must satisfy
the {\sl first binomial poset compatibility\/} condition
 $$1 = a_1 \le a_2 \le a_3 \le \ldots, \qquad\hbox{and}\qquad
{\prod\3l{j+1}{j+i}a_l\over\prod\3l1ia_l} \in \N = \{0,1,2,\ldots\}
\leqno(1)$$ for all $i,j \in \P = \{1,2,3,\ldots\}$.

 \smallskip
 Note that by the Ehrenborg-Readdy theorem
\cite{Ehrenborg_Readdy_Eulerian}, there are but two
isomorphism classes of Eulerian (strongly confluent infinite
binomial) posets, with
the atomic sequences $(1,2,2,2,2\ldots) = (1,2^\infty)$ and
$(1,2,3,4,5,\ldots) = (i)_i$, respectively. In contrast,
we may classify the strongly confluent infinite binomial posets with $\A =
(1^2,2^\infty) = (1,1,2,2,2,\ldots)$, by means of the infinite binary
strings without consecutive 1's. In particular, there is an
uncountable number of isomorphism classes of them.

 \medskip
 Formally, given a poset $X$ as above, we shall define a
string $\phi(X) = (\iota_j)\3j1\infty \in [2]^\P$, such that
$\phi(X) = \phi(X') \iff X \simeq X'$, and that $S := \mathop{\rm Im}
\phi = \bigl\{ (\iota_j)_j \in [2]^\P : \max\limits_{j\in\P}
(\iota_j+\iota_{j+1}) \le 3\bigr\}$.

 \section{The main results}
 For any (strongly confluent infinite downwards bounded) binomial
poset $X$, let $X_i := \{x\in X: \rank(x) = \length(\2\0x) = i\}$.
 For $0 \le i \le j$, let
the \emph{section} $X_{i,j}$ be $\bigcup\limits\3lijX_l$. Name $X$
\emph{of bounded atomic number type}, if $a := \lim \A :=
\lim_{n\rightarrow\infty} a_i < \infty$. Call \A\ {\sl realised\/} by
$X$. We have the following lemma:

\proclaim Lemma. For a strongly confluent infinite downwards bounded poset
$X$ with atomic sequence \A, et cetera, as above, we have
\clause{$(a)$}{The following conditions are equivalent:\br
 $i$. $a<\infty$;\br
 $ii$. $\ca X_i < \infty$ for some $i\in\P$;\br
 $iii$. All $\ca X_i < \infty$;\br
 $iv$. $\sup_i \ca X_i < \infty$.}
 \clause{$(b)$}{If indeed $X$ is of bounded atomic number type (i.e.,
 the conditions in $(a)$ are fulfilled), then
 $$\ca X_i = {a^i\over B(i)};\quad \sup_i \ca X_i =
\lim_{i\rightarrow\infty} \ca X_i = \prod\3i1\infty (aa_i^{-1})\,;$$
$X$ is countable; and every element in $X$ is covered by exactly $a$
elements.}

 {\sl Proof}. By the confluency condition, there is a chain $(x_i)_i$,
such that for any finite sequence
of different elements $\row yr \in X_1$, there is an $i_r \in \P$, such
that \row yr cover \0 in \2\0{x_{i_r}}; and without loss of
generality, we may assume that $x_{i_r} \in X_{j_r}$ for some $j_r \ge
i_r$. Hence, if there is an infinite sequence $y_1,y_2,\ldots \in X_1$,
then on the one hand $a = \lim\limits_i a_{j_i} \ge \lim\limits_i i =
\infty$, while on the other $\ca X_i \ge \ca X_1 $ and thus is
infinite, for each $i\in \P$.

 Conversely, assume that $X_1 = \{\row yr\}$.
Then $a_i = \ca X_1 $ for
$i \ge j_r$, whence $a = \ca X_1 $, indeed. Moreover, similar
arguments hold for the elements covering any fixed $x \in X$,
instead of $X_1$; whence then indeed $x$ is covered by $a$ elements.
In particular, there are exactly $a^i$ saturated chains of length
$i$, starting at \0. On the other hand, such a chain is a maximal
chain in \2\0x for some $x \in X_i$; whence there are $\ca X_i B(i)$
different such chains. In particular indeed $\ca X_i = a^iB(i)^{-1}
< \infty$. Now, the remaining claims follow easily.\square

 \smallskip
 Now, fix $a_1=a_2=1$, $a_3=~\ldots~=a = 2$. Then $\ca X_1 = 2$,
and $\ca X_2 = \ca X_3 =\ldots=4$. For any $i \in \P$, (the Hasse
graph of) the section $X_{i+1,i+2}$ is a bipartite, 2-regular
graph on 8 vertices, and thus is isomorphic to either \c8 or to
$\c4\cup\c4$ (distinguishable by the number of components). Put
$\iota_i :=  3 - \comp X_{i+1,i+2}$, and
$\phi(X) := (\iota_i)_i$.

 \smallskip
 Obviously, $\bigl(X \simeq X' \implies \phi(X) = \phi(X')\bigr)$.

 \smallskip
 For the moment, fix $i \in \P$.\br
 Recall that there are three possible $2+2$ partitions of a given
4-set. Now, any $y \in X_{i+2}$ covers exactly two elements in
$x_{i+1}$; and if $y$ covers $x_1$ and $x_2$, then there is a $y' \in
X_{i+2}$ that covers $X \setminus \{x_1,x_2\}$, defining a $2+2$
`cover' partition $\bigl\{\{x_1,x_2\},\;X\setminus\{x_1,x_2\}\bigr\}$
of $X_{i+1}$. Thus, $X_{i+1,i+2}$ induces $3-\comp X_{i+1,i+2}$
different cover partitions of $X_{i+1}$; and symmetrically equally
many `co-cover' partitions of $X_{i+2}$. Likewise, $X_{1,2}$ induces
one co-cover partition of $X_2$. More precisely, up to isomorphisms
$X_{0,2} = \{\0,z_1,z_2,w_1,w_2,w_3,w_4\}$ with the cover relations
$\0<z_\nu$, and $z_\nu<w_\mu$ for $\mu \equiv \nu \pmod2$; inducing
the partition $\bigl\{\{w_1,w_3\},\;\{w_2,w_4\}\bigr\}$. Given this,
and given all {\sl sets\/} $X_j$ (and again letting $i$ float), we get

 \proclaim Lemma. A definition of \c8 or $\c4\cup\c4$ covering
structures for all segments $X_{i+1,i+2}$ ($i\in\P$) defines a
type $(1,1,2,2,2,\ldots)$ confluent binomial poset structure on $X$,
iff there is no $i$ and $2+2$ partition of $X_{i+1}$ which is at the
same time an induced cover and an induced co-cover partition.

 {\sl Proof}. If the partition $\bigl\{\{b,c\},\; \{d,e\}\bigr\}$ of
$X_{i+2}$ is both a cover and a co-cover partition, then there are
$f \in X_i$ and $g\in X_{i+2}$, such that the interval $[f,g] =
\{f,b,c,g\}$, contradicting $A(2) = 1 \ne 2$. Thus the partition
avoidance condition is necessary. For the sufficiency, inspect all
interval maximal chain numbers, employing the fact that the partition
avoidance condition in particular implies that\br $\bigl( (j-i \ge 2$
and $b \in X_i$ and $c \in X_j) \implies b<c \bigr)$.\square

 \smallskip
 With $X$ and $(\iota_i)_i$ as above, the lemma immediately yields
 that $\iota_i+\iota_{i+1} \le \ca \{2+2\;\hbox{partitions of
}X_{i+2}\} \allowbreak\le 3$; but that this condition is sufficient
for the existence of $X$. Finally, the implication $(\phi(X) =
\phi(X') \implies X \simeq X')$ may be proved by constructing
isomorphisms $X_{0,i} \simeq X'_{0,i}$ by induction on $i$. Thus
the claim is proved.

 \proc Remark. Note, that the two kinds of sections, \c8 and
$\c4\cup\c4$ were employed in order to construct non-isomorphic
strongly confluent finite binomial posets with the same `factorial
functions', or with non-isomorphic intervals of the same size,
already in \cite[Section 8.2]{D_R_S}; cf.\ e.g.\ Figure~5 (p.~309)
therein.

 \medskip
 The following figure illustrates the posets with $\phi = (1^\infty)$
or $\phi = (1,2,1,2,1,2,\ldots)$, respectively. In a sense, these
are the extremal possibilities. There are as many non-isomorphic
posets (with atomic sequence $\A = (1^2,2^\infty)$) as there are
different $\phi$, i.e., $\aleph_0$ many. Moreover, there is a $\phi$,
say $\overline\phi$, with the `versal property' to contain each
finite substring of $[2]^\infty$ without consecutive 2's. (You may
e.g.\ enumerate all such substrings as $s_1,s_2\kdots$ and put
$\overline\phi := (s_1,1,s_2,1,s_3,\ldots)$.) The corresponding
poset contains copies of all possible intervals in
posets with the given \A; and in particular contains a Fibonacci
number of non-isomorphic intervals of any given positive length.

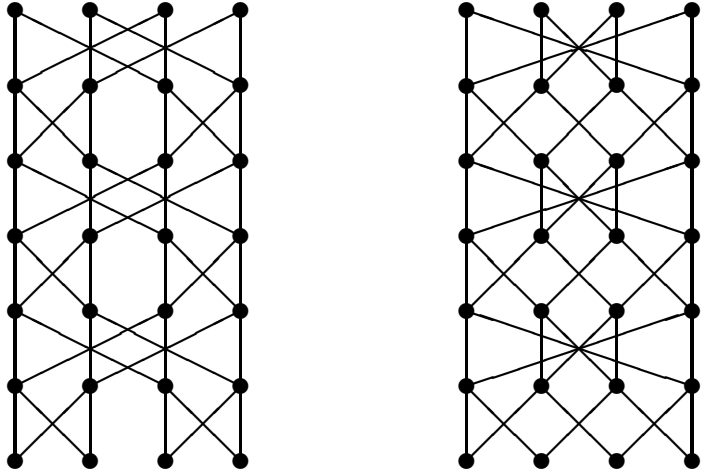
\setlength{\unitlength}{1.00mm}
\begin{picture}(0,0)(90,70)

\thicklines
\put(120,0){
\put(0,0){
\multiput(0,0)(10,0){4}{\multiput(0,0)(0,10){7}{\circle*{2}}}

\multiput(0,0)(0,20){3}{
\put(0,0){\line(0,1){10}}
\put(10,0){\line(0,1){10}}
\put(0,0){\line(1,1){10}}
\put(10,0){\line(-1,1){10}}
\put(20,0){\line(0,1){10}}
\put(30,0){\line(0,1){10}}
\put(20,0){\line(1,1){10}}
\put(30,0){\line(-1,1){10}}
}

\multiput(0,10)(0,20){3}{
\put(0,0){\line(0,1){10}}
\put(20,0){\line(0,1){10}}
\put(0,0){\line(2,1){20}}
\put(20,0){\line(-2,1){20}}
\put(10,0){\line(0,1){10}}
\put(30,0){\line(0,1){10}}
\put(10,0){\line(2,1){20}}
\put(30,0){\line(-2,1){20}}
}
}

\put(60,0){
\multiput(0,0)(10,0){4}{\multiput(0,0)(0,10){7}{\circle*{2}}}

\multiput(0,0)(0,20){3}{
\put(0,0){\line(0,1){10}}
\put(0,0){\line(1,1){10}}
\put(10,0){\line(-1,1){10}}
\put(10,0){\line(1,1){10}}
\put(20,0){\line(-1,1){10}}
\put(20,0){\line(1,1){10}}
\put(30,0){\line(-1,1){10}}
\put(30,0){\line(0,1){10}}
}

\multiput(0,10)(0,20){3}{
\put(0,0){\line(0,1){10}}
\put(0,0){\line(3,1){30}}
\put(30,0){\line(-3,1){30}}
\put(30,0){\line(0,1){10}}
\put(10,0){\line(0,1){10}}
\put(10,0){\line(1,1){10}}
\put(20,0){\line(-1,1){10}}
\put(20,0){\line(0,1){10}}
}
}
}
\end{picture}

\vspace*{70mm}

 \section{Further questions}

 This section contains some comments and questions that the
`experts' couldn't answer right on the spot, but seem to have
thought of. Thus, I strongly suspect that the results I
mention, or similar examples, are known (even if not published).

 \smallskip
  As is well-known, the first binomial
compatibility  condition (1) is not sufficient to guarantee
realisability. I shall provide a brief proof, based on the
following two lemmata:

 \proclaim Lemma. The infinite sequence \A\ is realisable iff it is
realised by some strongly confluent binomial poset.

 {\sl Proof}. Suppose \A\ is realised by $P$. For any $n \in \N$,
let $I_n$ be the set of isomorphism classes of intervals of
length $n$ in $P$. Each such interval has a fixed finite size,
whence each $I_n$ is finite. For each $n > 0$ and $Q \in I_n$,
choose one element $f(Q) \in I_{n-1}$, such that any interval \2xy
of type $Q$ has a subinterval \2xz of type $f(Q)$. Make $I :=
\bigcup_nI_n$ to an infinite directed tree, by letting the children
of $Q \in I_n$ be $f^{-1}(Q) \subseteq I_{n+1}$. By K\"onig's lemma,
there is an infinite path in $I$. The direct limit of this path
(in the natural sense) indeed is
a strongly confluent binomial poset realising \A.\square

 \proclaim Lemma. If the finite sequence $(a_1, a_2\kdots a_n)$
satisfies (1) for all $i,j\in\P$ such that $i+j\le n$, then it may be
extended to an infinite sequence \A\ satisfying (1); and \A\ may
be chosen with $\lim \A < \infty$.

 \proc Proof. Put $a_{n+1} := a_{n+2} := \ldots :=
\mathop{\rm lcm} (\row an)$.\square

 \smallskip
 Now, note that $(1,2,3,4,4)$ fulfils the assumptions of the lemma,
but cannot be the atomic numbers sequence for any binomial poset
interval $X = [\0,\1]$ of rank~5. In fact, for any integers $n \ge 3$
 and $a_n \ge n-1$, the sequence $(1,2\kdots n-1,a_n)$ is realisable
if and only if $n$ divides $a_n$. (If $a_n = kn$, say, with $k \in
\P$, then we may construct such an interval by stripping the boolean
lattice on $n$ atoms of its top and bottom elements, taking $k$
disjoint copies of the resulting poset, and adding new top and bottom
elements. Conversely, if we have a realising interval $X =
\2\0\1$ of the sequence, then we may define a binary relation $R$ on
the $a_n$-set $X_1$ by $xRx' \iff \exists\,y \in X_2 : x,x' \in \2\0y$.
Now, by means of the fact that any proper subinterval of $X$ is
boolean,, we may prove that $R$ is an equivalence relation, and that
each equivalence class has the size $n$. Thus, if $k$ is the number
of equivalence classes, indeed we get $a_n = kn$.)

 Thus and by the lemmata, there is an infinite sequence $\A = (1, 2,
3, 4, 4,\ldots)$ and fulfilling (1) but non-realisable as the
atomic numbers sequence of any binomial poset whatsoever. In fact,
the sequence $(1,2,3,4^2,6^\infty)$ satisfies (1).

 Another example is provided by the fact that the atom number
sequence $(a_1,a_2,a_3) = (1,m,m+1)$ is (uniquely) realised by the
interval $\{\0;\row x{m+1};\row y{m+1};\1\}$, with $x_i<y_j$ if and
only if $i \ne j$; but that for $m \ge 3$ this is not extendable to
any binomial interval of larger length. Thus, $\bigl( 1,m,m+1,
\bigl(m(m+1)\bigr)^\infty \bigr)$ satisfies (1), but is not
realisable.

 \medskip
 It  would be interesting to find further general conditions that
(finite or infinite) atomic number sequences must fulfil in order to
be realisable; and optimally a necessary and sufficient set of
numerical conditions.

 \smallskip
 It is fairly easy to see that $(1) \;\&\; \lim \A \le 3$ imply that \A\
is realisable. In fact, any prime number \A\ limit is covered by
 \proclaim Lemma. For any $m,n\in\P$, the sequence $(1^m,n^\infty)$
is realised by $X = \bigcup_i X_i$, defined by $X_i :=
\{i\} \times [n]^{\min(i,m)}$, with a cover relation $(i,s) < (i+1,t)$
iff the maximal proper left substring of $t$ is a right substring of
$s$.\square

 It is also fairly easy to
extend this construction to any sequence $\A = (a_i)\3i1\infty \in
\{1\}\times\P^\infty$ such that $a_i|a_{i+1}$ for all $i\in\P$. In
particular, this covers any sequence satisfying (1) and with a
prime power limit. Thus and by the example {\sl supra}, the smallest
limit $a$ for which (1) is not a sufficient condition is $a = 6$. I
suspect that there are similar counterexamples for any limit $a$
which has different prime factors.


\section{Acknowledgements}

 I thank Margaret Readdy, Richard Stanley, and Richard Ehrenborg for
introducing me to this subject, in seminairs and stimulating conversation
at the Mittag-Leffler Institute; and especially the latter Richard for
\TeX nical assistance with the illustration.

\newcommand{\journal}[6]{{\sc #1,} #2, {\it #3} {\bf #4} (#5), #6.}
\newcommand{\book}[4]{{\sc #1,} ``#2,'' #3, #4.}
\newcommand{\preprint}[3]{{\sc #1,} #2, preprint #3.}
\newcommand{\JCTA}{J.\ Combin.\ Theory Ser.\ A}


 \bigskip
\noindent
 Matematiska Institutionen \\
   Stockholms Universitet \\
   SE-106~91 Stockholm\\ SWEDEN\\
{\tt joeb@math.su.se}

\end{document}